\documentclass[11pt]{article}
\usepackage{amssymb,latexsym,amsmath,amsbsy,amsthm,amsxtra,amsgen,dsfont,pdfsync,color,graphicx}
\usepackage{mathrsfs}
\newfont{\msbm}{msbm10 at 11pt}

\newfont{\msbmsm}{msbm10 at 8pt}

\newtheorem{Theo}{Theorem}[section]

\newtheorem{corollary}[Theo]{Corollary}

\newtheorem{Exm}[Theo]{Example}

\newtheorem{remark}[Theo]{Remark}
\newtheorem{conjecture}[Theo]{Conjecture}

\begin{document}
\title{Probabilistic aspects of $\Lambda$-coalescents \\ in equilibrium  and in evolution}

\author{G\"otz Kersting\thanks{Institut f\"ur Mathematik, Goethe Universit\"at, Frankfurt am Main, Germany \newline kersting@math.uni-frankfurt.de, wakolbinger@math.uni-frankfurt.de \newline Work partially supported by the DFG Priority Programme SPP 1590 ``Probabilistic Structures in Evolution''}$\ $ and Anton Wakolbinger$^*$}
\date{}
\maketitle

\begin{abstract}
We present approximation methods which lead to law of large numbers and fluctuation results for functionals of $\Lambda$-coalescents, both in the dust-free case and in the case with a dust component. Our focus is on the tree length (or total branch length) and the total external branch length, as well as the  time to the most recent common ancestor and the size of the  last merger. In the second part we discuss evolving coalescents. For certain Beta-coalescents we analyse  fluctuations of a class of functionals in appropriate time scales. Finally we review results of Gufler on the representation of evolving $\Lambda$-coalescents in terms of the lookdown space.
\end{abstract}

\section[Introduction]{Introduction}

$\Lambda$-coalescents \index{Lambda@$\Lambda$-coalescent}\index{coalescent!Lambda@$\Lambda$-}are a well-established model for the genealogy of a sample of individuals out of a large population. We refer to Birkner and Blath \cite{GKAW-BB20} for a comprehensive introduction into the biological background, and a discussion of statistical aspects not only of $\Lambda$-coalescents but also of \index{Xi@$\Xi$-coalescent}\index{coalescent!Xi@$\Xi$-} $\Xi$-coalescents. The latter are \index{multiple merger coalescent}\index{coalescent!multiple merger} mutiple merger coalescents which, unlike the $\Lambda$-coalescents, also admit {\em simultaneous} multiple \index{merger} mergers, and constitute the most general class of exchangeable coalescents. F. Freund \cite{GKAW-Fr20} emphasises various aspects of $\Xi$-coalescents, and A. Sturm \cite{GKAW-St20} treats diploid models in connection with $\Xi$-coalescents.

$\Lambda$-coalescents are Markov processes on the set of partitions of the natural numbers $\mathbb N$, and their law is determined by a finite measure $\Lambda$ on the interval $[0,1]$. 
Their dynamics is readily described by  their restrictions on the finite sets $[n]:=\{1,2,\ldots,n\}$, $n\in \mathbb N$.  These projections, which are called \index{coalescent!n@$n$-}  {\em$n$-coalescents}, are Markov processes on the finite set of partitions of $[n]$ with the initial value $\{\{1\}, \{2\}, \ldots,\{n\}\}$ and with the property that  $k$ specified blocks out of a given partition $\{B_1, \ldots, B_b\}$ of $[n]$   merge at rate
\[ \lambda_{b,k} := \int_{[0,1]} p^k(1-p)^{b-k} \, \frac{\Lambda(dp)}{p^2} \ , \ 2 \le k \le b . \]
There is an equivalence between the fact that there exists a finite measure $\Lambda$ such that the last equations holds and the fact that the process is consistent in the sense that the restriction to $[n]$ of the $(n + 1)$-coalescent is the $n$-coalescent.

In biological terms we may regard an $n$-coalescent as a tree with $n$ leaves and a root, representing  the genealogy of the sample of  $n$ individuals (or genes) and their most recent common ancestor, respectively. Among the tree's branches one distinguishes between the \index{branch!total external}$n$   {\em external ones}, each of them connecting a leaf to its neighbouring  internal node, and the \index{branch!total internal}{\em internal branches}.  The  branches have lengths which represent  time durations in the genealogy. 

Many special cases of $\Lambda$-coalescent are important and have been studied in greater detail, such as the \index{Kingman coalescent}\index{coalescent!Kingman}Kingman coalescent, the \index{Bolthausen-Sznitman coalescent}\index{coalescent!Bolthausen-Sznitman}Bolthausen-Sznitman coalescent or the more comprehensive class of \index{Beta-coalescent}\index{coalescent!Beta-}Beta-coalescents.
Motivated by biological applications a number of functionals of such coalescents have been analysed in the literature.  Our results will concern  general $\Lambda$-coalescents, with a focus on the trees'  \index{branch length!total}\index{tree length}total and \index{branch length!total external}external branch lengths.   We also present results on  the time to the \index{most recent common ancestor}most recent common ancestor  (MRCA) and the size of the \index{merger!size of last}\index{coalescent!last merger} coalescent's last merger.

\begin{figure}[t]
    \begin{center}
   \includegraphics[height=4cm]{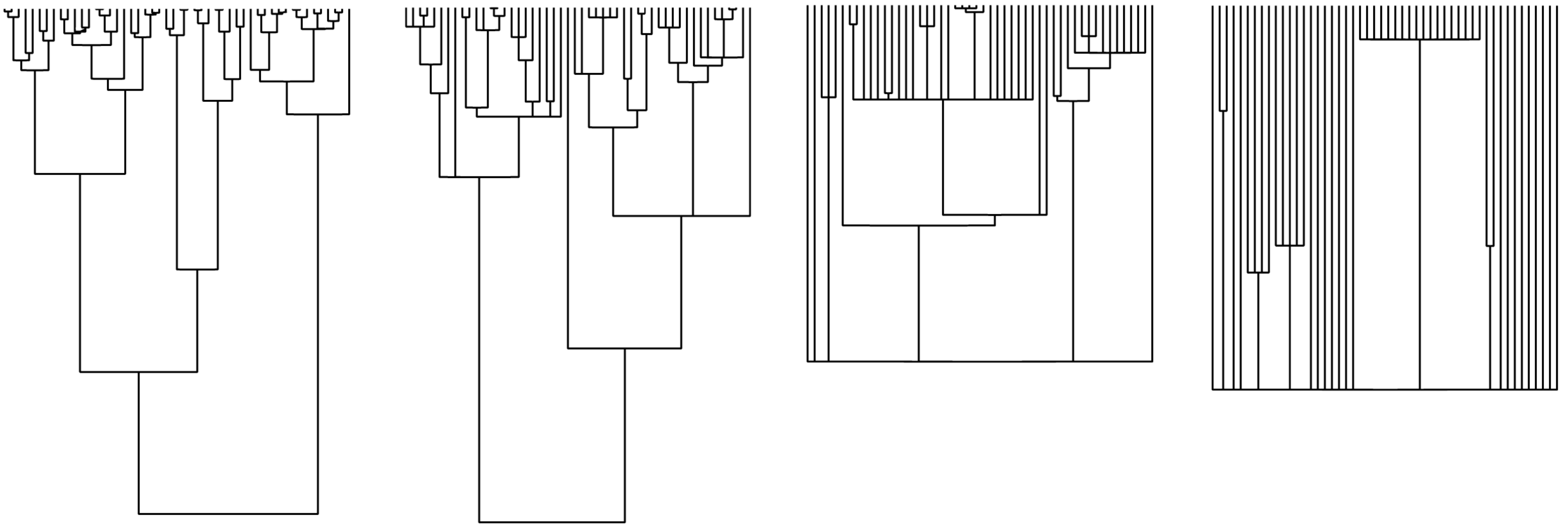}
\end{center}
\caption{\label{GKAW-fig0} Four realisations of $n$-coalescents with $n=50$ from the family of \index{Beta-coalescent}\index{coalescent!Beta-}Beta($2-\alpha, \alpha$)-coalescents. From left to right: $\alpha=2$ \index{Kingman coalescent}\index{coalescent!Kingman}(Kingman),  $\alpha=3/2$, $\alpha=1$ \index{Bolthausen-Sznitman coalescent}\index{coalescent!Bolthausen-Sznitman}(Bolthausen-Sznitman), $\alpha=1/2$. The right-hand tree stems from a coalescent with a dust component,  and  the two left-hand trees belong to coalescents coming down from infinity. }
\end{figure}
As illustrated in Figure \ref{GKAW-fig0}, coalescent trees can have very different behaviours depending on the choice of the measure $\Lambda$.
The key quantity in this context turns out to be the \index{block-counting process!rate of decrease}coalescent's  {\em rate of decrease} given by the rates
\[ \mu(b) := \sum_{k=2}^b (k-1) \binom bk \lambda_{b,k}\ , \ b \ge 2 . \]
This formula reflects the fact that the coalescence of $k$ blocks results in a decrease of $k-1$ in the number of blocks. Some of the most fundamental properties of $\Lambda$-coalescents can
be characterised in terms of the $\mu(b), b \ge 2$. A $\Lambda$-coalescent is said \index{dust component}\index{dust component!lack of}to {\em lack a dust component} (or  simply to \index{dust-free}be {\em dust-free}) if the probability   that  there are external branches within the $n$-coalescents  reaching right down to the tree's root is asymptotically vanishing in the limit $n \to \infty$. This feature takes place if and only if the \index{block-counting process!rate of decrease per capita}{\em rate of decrease per capita} is divergent, i.e.  
\[ \frac{\mu(b)}b \to \infty \]
as $b \to \infty$. Originally this kind of behaviour has been characterized by the condition $\int p^{-1} \, \Lambda(dp)=\infty$ \cite{GKAW-P99}, its equivalence to the above requirement is shown in \mbox{\cite[Lemma 1 (iii)]{GKAW-DiKe18a}}. 
Another property of importance is the stronger notion that a $\Lambda$-coalescent \index{coming down from infinity}{\em comes down from infinity}, which means that for any $t>0$ the number of lineages at time $t$ in the $\Lambda$-coalescent is finite a.s. Expressed in terms of $n$-coalescents this means that for any $t>0$ the numbers of their lineages at time $t$ are tight as $n \to \infty$. This feature arises \cite{GKAW-S00} if and only if
\[ \sum_{b=2}^\infty \frac 1{\mu(b)} < \infty.\]

The paper is subdivided in two parts. The first part, consisting of Sections~\ref{Nodust}, \ref{Dust} and \ref{asexBeta}, focusses on aspects of the genealogy of $n$ individuals that live at some fixed time, and on their description by coalescent functionals. Sections \ref{Nodust} and \ref{Dust} treat general $\Lambda$-colaescents. Each of them presents a different approximation method, one for a class of functionals of the $\Lambda$-coalescent and the other for the block-counting processes. In Section \ref{asexBeta} we turn to Beta-coalescents, and present an asymptotic expansion for a further class of functionals.

In the second part we are  going to describe and analyse \index{process!coalescent-valued}\index{coalescent-valued process}{\em coalescent-valued processes} that are embedded in evolving populations. Here we distinguish two different perspectives. In Section \ref{evolvn}  we think of $n$ as a population size that is constant over time. Then each of the functionals studied in the first part gives rise to a real-valued process, and we are going to obtain limit laws for  sequences of such processes as $n \to \infty$. There, the factor that translates \index{time scale!generation}the {\em generation time scale} $s$ into the \index{time scale!evolutionary}{\em evolutionary time scale} $t$ will vary with $n$, cf. Remark \ref{remarkgentime}.
In Section~\ref{EvolvLambda} we will consider genealogies of infinite populations in their evolutionary time scale. As we shall see, these genealogies can be viewed as \index{evolving coalescent}\index{coalescent!evolving}\index{coalescent!i@$\infty$-}evolving $\infty$-coalescents. Here $n$ looses its meaning as a total population size but still can be seen as the size of a sample taken from the population at some fixed time. 

\section[Dust-free $\Lambda$-coalescents]{An approximation method for \index{dust-free}dust-free $\Lambda$-coalescents}\label{Nodust}

\subsection[A class of functionals of $\Lambda$-coalescents]{A class of functionals of $\Lambda$-coalescents}
We start by introducing some notation. For fixed $n\in \mathbb N$, let $N_n=(N_n(u))_{u \ge 0}$ denote the \index{block-counting process}block-counting process. Thus $N_n(u)$ gives the number of blocks within the $n$-coalescent at time $u$, in particular $N_n(0)=n$. The process $N_n$ is a Markov pure-death process with state space $[n]$, with an absorbing state 1 and \index{block-counting process!jump rate}with a jump rate in state $b$ given by
\[ \lambda(b) := \sum_{k=2}^b  \binom bk \lambda_{b,k}\ , \ b \ge 2 . \]
The   Markov chain $X = X^{(n)}$ is the decreasing path $n=X_0^{(n)}> X_1^{(n)} > \cdots > X_{\tau_n}^{(n)} = 1$, embedded into the Markov process $N_n$, with $\tau_n$ being the total number of merging events. The waiting time in the state $X_i=X_i^{(n)}$, $i=0,1,\ldots$, is denoted by $W_i=W_i^{(n)}$. For the sake of readability we will suppress the superscript $n$ in   $X_i$ and $W_i$.

A number of functionals of $n$-coalescents can be expressed or closely approximated by quantities of the form
\begin{align} F_n(X):=\sum_{i=0}^{\rho_n-1} f(X_i) 
\label{Fn}
\end{align}
with some function $f:[2,\infty)\to \mathbb R$ and some random variable $0 \le \rho_n \le \tau_n$. Here are a few examples.

\begin{Exm} {\bf Number of coalescences.} If we \index{coalescent!number of coalescences}set $f \equiv 1$ and $\rho_n=\tau_n$, then $F_n(X)$ equals the number $\tau_n$ of merging events.
\end{Exm}
\begin{Exm}\label{exabs}{\bf Absorption time.} For the \index{coalescent!absorption time}\index{absorption time}time  to the most recent common ancestor, $\tilde \tau_n:=W_0+ \cdots+ W_{\tau_n-1}$, we have
\[ \mathbb E[ \tilde \tau_n \mid X] = \sum_{i=0}^{\tau_n-1} \frac 1{\lambda(X_i)} \ . \]
Here the choice  $f(x):=\lambda(x)^{-1}$  in \eqref{Fn} leads to  an effective approximation of~$\tilde\tau_n$.
\end{Exm}
\begin{Exm}{\bf Tree length.} \index{tree length}\index{branch length!total}Similarly the total length \[\ell_n= X_0W_0+ X_1W_1+ \cdots + X_{\tau_n-1}W_{\tau_n-1}\] of all branches can be well approximated by
\begin{align} \mathbb E[\ell_n \mid X] =  \sum_{i=0}^{\tau_n-1} \frac {X_i}{\lambda(X_i)} \ . 
\label{treelength}
\end{align}
\end{Exm}
\begin{Exm} {\bf External branches.} \label{exam} \index{branch!total external}The choice $f(x):= 1/x$, $x \ge 2$, deserves particular interest. Let $\zeta^n$ be the number of  coalescing events before an external branch chosen at random out of the $n$ possible
merges with some other branches within an $n$-coalescent. We have  the fundamental formula
\[ \mathbb P(\zeta_n \ge k \mid X) =\frac{X_k-1}{n-1} \prod_{i=0}^{k-1} \Big(1- \frac 1{X_i}\Big) \ \text{ a.s.} \]
for any $k \in \mathbb N$ (see \cite[Lemma 4]{GKAW-DiKe18a}). This implies 
\begin{align} \label{zetan}
\mathbb P( \zeta_n\ge k \mid X)  = \frac{X_k-1}{n-1} \exp \Big(-\sum_{i=0}^{k-1} \frac 1{X_i}+ O(X_{k}^{-1})\Big)\  \text{ a.s.}  
\end{align}
Here  the functional $F_n(X)$ with $f(x) = 1/x$   is  located in the exponent. With this approximation  one gains access to the lengths of external branches, see \cite{GKAW-DiKe18a,GKAW-DiKe18b}.
\end{Exm}
\subsection{An approximation method for dust-free $\Lambda$-coalescents}

Following \cite{GKAW-DiKe18a} we present an approach to obtain laws of large numbers for the random variables $F_n(X)$ from \eqref{Fn}.  This approach relies on the following intuition. We note that the sequences $\lambda(b)$ and $\mu(b)$, $b \ge 2$ can be naturally extended to smooth functions $\lambda,\mu:[2,\infty) \to \mathbb R$, see \mbox{\cite[Equations (4) and (5)]{GKAW-DiKe18a}}. Let
\[ \Delta_i:=X_i-X_{i-1} \ , \quad  \nu(x):= \frac{\mu(x)}{\lambda(x)} \ . \]
Then we have the following approximation in two steps
\begin{align} F_n(X) \approx \sum_{i=0}^{\rho_n-1} f(X_i)\frac{\Delta_{i+1}}{\nu(X_i)} \approx \int_{r_n}^n f(x) \frac{dx}{\nu(x)} \ ,
\label{approx}
\end{align}
where $(r_n)_{n \ge 1}$ is a sequence of real numbers, $\rho_n:= \min\{ k \ge 0: X_k < r_n\}$ (and the  symbol $ \approx$ is understood in a heuristic manner). The right-hand part of \eqref{approx} can be seen as a Riemann approximation of the integral. To provide a sufficiently good fit, a natural requirement is that $ \sup_{\rho_n\le  i\le n} \tfrac{\Delta_{i+1}}{X_i}$ becomes small in probability as $n \to \infty$. Therefore, in order to avoid very large jumps, we confine ourselves to a small-time regime meaning that the time $\tilde \rho_n:= \inf \{ u \ge 0: N_n(u) < r_n\}$ converges to 0 in probability as $n \to \infty$. For coalescents \index{coming down from infinity}coming down from infinity this simply means that $(r_n)$ diverges. Otherwise this is a stronger restriction, then $(r_n)$ has to diverge sufficiently fast. 

The rationale for the left-hand approximation in \eqref{approx} rests in the observation that we have
\[ \mathbb E[\Delta_{i+1}\mid X_i] = \nu(X_i) \ \text{a.s.}\]
and consequently the difference of both left-hand terms may be embedded in the martingale $M=(M_k)_{k \ge 0}$ given by
\[ M_k := \sum_{i=0}^{k\wedge \tau_n -1}f(X_i) - \sum_{i=0}^{k\wedge \tau_n -1}f(X_i)\frac{\Delta_{i+1}}{\nu(X_i)} \ , \ k \ge 0 \ . \]
Its quadratic variation may be bounded in the regime $ \sup_{\rho_n\le  i\le n} \tfrac{\Delta_{i+1}}{X_i}= o_P(1)$ by means of the estimate
\begin{align*} \sum_{i=0}^{\rho_n-1} f(X_i)^2\frac{\Delta_{i+1}^2}{\nu(X_i)^2}&= o_P\Big(\sum_{i=0}^{\rho_n-1} f(X_i)^2\frac{\Delta_{i+1}X_i}{\nu(X_i)^2}\Big)\\& = o_P\Big(\max\limits_{r_n \le x \le n } \frac{xf(x)}{\nu(x)}\sum_{i=0}^{\rho_n-1} f(X_i)\frac{\Delta_{i+1}}{\nu(X_i)}\Big)\ . 
\end{align*}
If  $(r_n)$  increases sufficiently slow then one can show 
\cite[Proposition 1]{GKAW-DiKe18a} that under suitable conditions on $f$  the maximum on the right-hand side of the previous display, $\max_{r_n\le x\le n} \frac{xf(x)}{\nu(x)}$, is $O\left(\int_{r_n}^n f(x)\frac{dx}{\nu(x)}\right)$,  and we end up with the estimate
\[ \sum_{i=0}^{\rho_n-1} f(X_i)^2\frac{\Delta_{i+1}^2}{\nu(X_i)^2} = o_P\Big(\Big(\int_{r_n}^n f(x) \frac{dx}{\nu(x)}\Big)^2 \Big), \]
which allows corresponding second moment estimates of the martingale $M$.

We note that martingale techniques have been applied earlier in \cite{GKAW-BBL10} to $\Lambda$-coalescents \index{coming down from infinity} 
coming down from infinity. There the speed of coming down from infinity was obtained by comparing the block-counting process (in continuous time) with the solution of an explosive ODE run backwards in time. 

The above described method does not cover this case, it follows a different strategy.
This kind of approach turns out to be applicable in great generality and appears to be a promising tool  for other decreasing Markov chains, too. For \mbox{$\Lambda$-coalescents} it leads to a law of large  numbers for the tree length and the total external branch length. To state this compactly, we use the following notation: For two sequences $(A_n)_{n \ge 1}$ and $(B_n)_{n \ge 1}$ of positive random variables we write $A_n \stackrel P \sim B_n$ and $A_n \stackrel 1 \sim B_n$, if the quotients $A_n/B_n$ converge to 1 in probability and in the $L^1$-norm, respectively.

\begin{Theo} \cite[Theorems 1 and 2]{GKAW-DiKe18a} Let $\ell_n$ be the \index{branch length!total}\index{tree length}total length of an $n$-coalescent and let $\bar \ell_n$ be the \index{branch length!total external}total length of its external branches. Then for a \index{dust-free}dust-free $\Lambda$-coalescent we have as $n \to \infty$
\[ \ell_n \stackrel 1 \sim  \int_2^n \frac x{\mu(x)} \, dx \quad \text{and} \quad \bar \ell_n \stackrel P \sim \frac {n^2}{\mu(n)} \ .\]
\end{Theo}
\mbox{}\\
The first part of the theorem was proven in \cite{GKAW-BBL14} for coalescents coming down from infinity, and was conjectured to hold for the larger class of dust-free $\Lambda$-coalescents.
 
It is natural to expect that an analogous theorem holds for \index{branch length!total internal}the {\em total internal branch lengths} $\hat \ell_n:= \ell_n-\bar \ell_n$ of  $n$-coalescents.  This is proven for a class of colalescents containing the \index{Bolthausen-Sznitman coalescent}\index{coalescent!Bolthausen-Sznitman}Bolthausen-Sznitman coalescent, see \cite[Theorem 3]{GKAW-DiKe18a}. Since there is further evidence that such a result holds in large generality,  we formulate the following conjecture.

\begin{conjecture} For a \index{dust-free}dust-free coalescent we have as $n \to \infty$
\[ \hat \ell_n \stackrel P \sim \int_2^n \Big(\frac x{\mu(x)}-\frac n{\mu(n)} \Big) \, dx \ . \]
\end{conjecture} 
\mbox{}\\
As indicated in Example \ref{exam}, our methodology is useful also to analyse single external branches in an $n$-coalescent. Here we have the following result.

\begin{Theo}\cite[Theorems 1.1 and 1.3]{GKAW-DiKe18b} Let $T_n$ be the length of an external branch chosen at random from an \index{coalescent!n@$n$-}$n$-coalescent. Then for a $\Lambda$-coalescent  without a \index{dust component} 
dust component  we have for all $u \ge 0$
\[   e^{-2u} +o(1)\le \mathbf P\Big( \frac {\mu(n)}n \, T_n \ge u \Big) \le \frac 1{1+u} +o(1) \]
as $n \to \infty$. 

Moreover, $\mu(n) T_n/n$ converges in distribution to a probability measure $\pi \neq \delta_0$ as $n \to \infty$, if and only if $\mu$ is a function varying regularly at infinity. Then its exponent $\alpha$ of regular variation fulfils $1 \le \alpha \le 2$ and we have for $\alpha =1$
\[ \pi(du) = e^{-u}\, du \]
and for $1<\alpha \le 2$
\[ \pi(du) =  \frac \alpha{(1+(\alpha-1)u)^{1+ \frac\alpha {\alpha-1}}}  \, du \ .\]
\end{Theo}
\noindent
The proof is based on formula \eqref{zetan}.
\begin{remark}\label{remarkgentime}
The second part of this theorem includes  results for the \index{Kingman coalescent}\index{coalescent!Kingman}Kingman and \index{Bolthausen-Sznitman coalescent}\index{coalescent!Bolthausen-Sznitman}the Bolthau\-sen-Sznitman coalescent  as to be found in the literature \cite{GKAW-W75,GKAW-DIMR07}. It suggests that $n/\mu(n)$ can be interpreted as the appropriate scaling of a generation's duration, i.e. the time at which a specific lineage out of the $n$ present ones takes part in a merging event, see \cite{GKAW-KWS14,GKAW-DiKe18b}.
\end{remark}

In addition to the previous theorem we point out that the lengths of the different external branches behave asymptotically  like i.i.d.~random variables, which for coalescents coming down from infinity \cite{GKAW-DiKe18b} also includes the external branches of maximal length.  For the  Bolthausen-Sznitman coalescent the picture changes:

\begin{Theo} \cite[Theorem 1.6]{GKAW-DiKe18b} Let $M_n$ denote the \index{branch length!longest external}length of the longest external branch in an $n$-coalescent. Then in case of a  \index{Bolthausen-Sznitman coalescent}\index{coalescent!Bolthausen-Sznitman}Bolthausen-Sznitman coalescent the random variables $\log \log (n)(M_n - t_n) $ converge in distribution 
as $n \to \infty$, with
\[ t_n := \log\log n - \log\log\log n + \frac {\log \log \log n}{\log\log n} \ . \]
The limit has the density 
$  ((1+e^u)(1+e^{-u}))^{-1}\ du$, $ u\in \mathbb R$, and is thus the standard logistic distribution.
\end{Theo}
\mbox{}\\
The latter arises as the distribution of the difference of two independent standard Gumbel random variables. Notably, if one proceeds to a point process description of the extremal lengths then it turns out that the limiting point process is a Poisson point process shifted by an independent standard Gumbel random variable. This random shift builds up away from the small time regime, consequently the analysis requires techniques different from the scheme \eqref{approx}. For details we refer to~\cite{GKAW-DiKe18b}.

As a final application of the above approximation methodology we present  the following result on the Bolthausen-Sznitman coalescent.

\begin{Theo} \cite[Theorem 4]{GKAW-DiKe18a}
Let $\bar \ell_{n,b}$ be the \index{branch length!branches of order $b$}total length of all branches of order $b\ge 1$ -- meaning the branches subtending  $b$ leaves {\em(}in particular $\bar \ell_{n,1}=\bar \ell_n${\em)}. Then we have for the  \index{Bolthausen-Sznitman coalescent}\index{coalescent!Bolthausen-Sznitman}Bolthausen-Sznitman coalescent  as $n\to \infty$
\[ \bar \ell_{n,1} \stackrel 1\sim \frac n{\log n}  \quad \text{and} \quad \bar \ell_{n,b} \stackrel 1\sim \frac 1{b(b-1)} \frac{n}{\log ^2 n} \]
for $b\ge 2$.
\end{Theo}
\noindent In the proof, formulae similar to \eqref{approx} come into play.
The result transfers immediately to the \index{site frequency spectrum}site frequency spectrum of the Bolthausen-Sznitman coalescent and is a counterpart to the  result in \cite{GKAW-BG08} on the corresponding allele frequency spectrum. For analoguous results for \index{Beta-coalescent}\index{coalescent!Beta-}Beta-coalescents  \index{coming down from infinity}coming down from infinity we refer to  \cite{GKAW-BBS07}.

\subsection{Fluctuations in $\Lambda$-coalescents: a conjecture}

The martingale approximation \eqref{approx} appears  to be fruitful also for the treatment of asymptotic fluctuations of the random variable $F_n(X)$ from \eqref{Fn}, in the following way:
\begin{align} F_n(X) \approx\int_{r_n}^n f(x)\frac{dx}{\nu(x)}+ \sum_{i=0}^{\rho_n-1} f(X_i) \frac{\nu(X_i)-\Delta_{i+1}}{\nu(X_i)}  \ .
\label{approx2}
\end{align}
This promises to improve results on the fluctuation behaviour, e.g.  for the  total length $\ell_n$ of $n$-coalescents, as follows.

Suppose that the $\Lambda$-coalescent is regularly varying with exponent $1 <\alpha < g$ with $g$ equal to the golden ratio $\frac 12 (\sqrt 5+1)$.  This means \cite{GKAW-DiKe18b} that
\begin{equation}\label{regvar}\int_{y}^1 \frac{\Lambda(dp)}{p^2} = y^{-\alpha}L(y^{-1}) \ , \quad 0<y<1 \ ,
\end{equation}
 where the function $L$ is slowly varying at infinity.  Then, using techniques from \cite{GKAW-DiKe18a}, one obtains that for any sequence $(a_n)_{n \ge 1}$ with $a_n \to \infty$ and $a_n=o(n)$ we have
\[ \mathbb P( \Delta_1 >a_n \mid X_0=n) = \frac{L(n/a_n)a_n^{-\alpha}}{\Gamma(2-\alpha)L(n)}(1+o(1))  \ , \] 
as $n \to \infty$, so we are in the domain of stable laws.
Choosing $(a_n)_{n \ge 1}$  such that
\begin{equation}\label{an}
\frac{L(n/a_n)a_n^{-\alpha}}{\Gamma(2-\alpha)L(n)}=\frac{1+o(1)}n\ , 
\end{equation}
the approximation \eqref{approx2} leads to  the
\begin{conjecture}{\bf Total length.} \index{tree length}\index{branch length!total}Assume \eqref{regvar} and \eqref{an}, with $1<\alpha<g$. Then
\[ \frac{\ell_n - \int_2^n \frac  x{\mu(x)} \, dx}{a_nn^{1-\alpha}/L(n)} \stackrel {\rm d}\longrightarrow -c  \zeta  \]
as $n\to \infty$, with
\[ c := \frac{(\alpha-1)^{1+1/\alpha}}{(1+\alpha-\alpha^2)^{1/\alpha}\Gamma(2-\alpha)}
 \]
and a stable random variable  $\zeta$  with index $\alpha$, which is  characterized by the properties $\mathbb E[\zeta]=0$, $\mathbb P(\zeta > z) \sim z^{-\alpha}$ and $\mathbb P(\zeta<-z)= o(z^{-\alpha})$ as $z \to \infty$.
\end{conjecture}
This conjecture is supported by corresponding results from \cite{GKAW-K12} on Beta-coalescents, where $L(n) \sim (\alpha \Gamma(\alpha)\Gamma(2-\alpha))^{-1}$ as $n \to \infty$.

\section[$\Lambda$-coalescents with a \index{dust component}dust component]{Approximating $\Lambda$-coalescents with a dust component}\label{Dust}

Following \cite{GKAW-KSW18} we  now come to a second approximation scheme, which is applicable beyond the small-time restriction addressed below formula \eqref{approx} and  extends thus to regions in the coalescent that are closer to the root. It is taylored to $\Lambda$-coalescents with a dust component, but has consequences  also in the dust-free case. 

As is well-known \cite{GKAW-P99} for $\Lambda$-coalescents with a dust component, the logarithm of the block-counting processes $N_n$ may be approximated by a \index{process!subordinator}\index{subordinator}subordinator $S$ in the sense that $\log N_n(t)$ and $\log n-S(t)$ are close to each other. Here the L\'evy measure of the subordinator is equal to the image of  the  measure $p^{-2}\Lambda(dp)$ under the mapping \mbox{$p \mapsto - \log (1-p)$.} The approximation relies on the fact that for this class of coalescents the very large mergers get dominant and smaller mergers may be neglected in the first instance. However, for various purposes this approximation is not good enough. The reason is that, if at a merging event $k$ out of $b$ blocks fuse, then the block-counting process does not decrease by the value $k$ but by $k-1$. Thus the \index{block-counting process}block-counting process decreases slightly slower than the accompanying \index{subordinator}\index{process!subordinator}subordinator. For the process $\log N_n$ this discrepancy has approximately size $1/N_n(u)$ at time $u>0$. Therefore we may state that, more accurately,
\[ \Delta \log N_n(u) \approx \log n-\Delta S(u) + \frac { \iota(N_n(u))\Delta u}{N_n(u)} \ , \]
where $\iota(b)$ denotes the subordinator's jump rate of those jumps which at state $N_n(u)=b$ find expression in the aforementioned discrepancy. This rate is given by $\int (1-(1-p)^{b} )p^{-2}\Lambda (dp)$. Thus, to make our intuition precise we introduce the process $Y_n$ as the solution of the stochastic integral equation
\[ Y_n(u)= \log n -S(u) + \int_0^u g\big(e^{Y_n(r)}\big) \, dr \]
with
\[ g(b):= \frac 1b\int_{[0,1]} (1-(1-p)^{b})\frac {\Lambda(dp)}{p^2} \ . \]
The integral takes finite values just for processes with a  dust component. 
Then we have the following theorem.

\begin{Theo}\label{dustapp}\cite[Theorem 10]{GKAW-KSW18}
If the $\Lambda$-coalescent has a \index{dust component}dust component then for all $\varepsilon >0$ there is an integer $k\ge 2$ such that for all $n$
\[ \mathbb P\big( \sup_{u \ge 0} |\log N_n(u)-Y_n(u)| I_{\{N_n(u) \ge k\}} \ge \varepsilon \big) \le \varepsilon \ .\]
\end{Theo}
\subsection{Time to absorption and size of the last merger}
Theorem~\ref{dustapp} is the essential tool for the proof of the next fundamental result, which holds for any $\Lambda$-coalescent.

\begin{Theo}\cite[Theorem 1]{GKAW-KW18}
The \index{absorption time}\index{coalescent!absorption time}absorption times $\tilde \tau_n$ defined in Example~\ref{exabs} satisfy \[ \tilde \tau_n \stackrel P\sim  \frac {\log n}\gamma \]
as $n \to \infty$, with
\[ \gamma := \int_{[0,1]} \big| \log  (1-p)\big| \frac{\Lambda(dp)}{p^2} \ .\]
\end{Theo}

\mbox{}\\
In the proof the case with dust is treated first, using Theorem \ref{dustapp}. Afterwards, dust-free coalescents (where $\gamma=\infty$) are approximated by coalescents with a \index{dust component}dust component. Along similar lines one also obtains a central limit theorem for $\tilde \tau_n$, see \cite{GKAW-KW18}.

Another application of this methodology concerns the \index{merger!size of last}size of the last merger of an \mbox{$n$-coalescents}. This is given by
\[ L_n:= X_{\tau_n-1}- 1\ . \]
For coalescents \index{coming down from infinity}coming down from infinity it is easy to see that this quantity converges in distribution to a limiting distribution on $\mathbb N$. The next theorems deal with the general case.

\begin{Theo}\cite[Theorems 1,2 and 3]{GKAW-KSW18}
For any $\Lambda$-coalescent fulfilling
\begin{align}
\int_{[0,1]}  |\log (1-p)|\, \Lambda(dp) < \infty 
\label{lastm}
\end{align}
 the sequence $(L_n)_{n\ge 1}$ is tight. Moreover,  in $\Lambda$-coalescents with a dust component the condition \eqref{lastm} is  necessary for tightness.

If additional to \eqref{lastm} we have  for all $d>0$ 
\begin{align} \Lambda\Big(\bigcup_{z=1}^\infty \{1-e^{-zd}\}\Big) < \Lambda((0,1]) \ ,
\label{Eldon}
\end{align}
then the sequence $(L_n)_{n\ge 1}$ is convergent in distribution.
\end{Theo}

\noindent
For \index{coalescent!Beta-}\index{Beta-coalescent}Beta-coalescents this result was obtained in \cite{GKAW-H15,GKAW-M14}. The counterexample \index{Eldon-Wakeley coalescent}\index{coalescent!Eldon-Wakeley}in \cite[Section 5: Non-convergence for   Eldon-Wakeley coalescents]{GKAW-KSW18} shows the relevance of condition \eqref{Eldon}. In case of convergence of the sequence $(L_n)_{n\ge 1}$  the coalescent \index{block-counting process!time-reversal}may be time-reversed in the following way. Let
\[ \hat N_n(u):= \begin{cases} N_n((\tilde \tau_n-u)-) & \text{for } u< \tilde \tau_n \ ,\\
n & \text{for } u \ge \tilde \tau_n \ , \end{cases} \]
in particular $\hat N_n(0)=L_n$. Again the next theorem covers also dust-free $\Lambda$-coalescents.

\begin{Theo} \cite[Theorem 5]{GKAW-KSW18}
If the sequence $(L_n)_{n \ge 1}$ converges in distribution, then the sequence of processes $(\hat N_n)_{n \ge 1}$ also converges in distribution in Skorohod space. The limiting process $\hat N_\infty$ is a Markov jump process with state space $\{2,3,\ldots\}$.
\end{Theo}

\noindent
For  formulas determining the jump rates of $\hat N_\infty$ we refer to \cite{GKAW-KSW18}. The special case of a \index{Bolthausen-Sznitman coalescent}\index{coalescent!Bolthausen-Sznitman}Bolthausen-Sznitman coalescent has been treated in \cite{GKAW-GM05}, and that of \index{Beta-coalescent}\index{coalescent!Beta-}Beta-coalescents in \cite{GKAW-H15}.

\section[An asymptotic expansion for Beta-coalescents]{An asymptotic expansion for Beta-coalescents}\label{asexBeta}
We now come to a different approximation method. It provides not only an approach to asymptotic distributions, but also  an asymptotic expansion in probability, which further allows   for the treatment of \index{coalescent!evolving}\index{evolving coalescent}{\em evolving} coalescents, see the next section.
We  specialize to the case in which $\Lambda$ is the \index{Beta-coalescent}\index{coalescent!Beta-}Beta$(2 - \alpha, \alpha)$-distribution for $1<\alpha < 2$. The following theorem gives an asymptotic expansion for a class of functionals, where the fluctuation term converges to a stochastic integral with respect to a compensated \index{process!Levy@L\'evy}\index{Levy@L\'evy process}L\'evy process that has an $\alpha$-stable distribution.  We will see that,  for the evolving Beta-coalescent, the same Poisson construction gives a representation of the fluctuation term as an $\alpha$-stable moving average process.

Let $\mathscr F$ be  the set of all differentiable functions $f$ which for some $c > 0$ and $0 < \zeta < \frac{1}{\alpha}$ obey 
$|f'(x)| \leq c x^{-\zeta - 1}$ 
for all $x \in (0, 1]$. For $n \in \mathbb N$ and $f\in \mathscr F$  let 
\begin{equation}\label{In}
\mathcal J_n(f):= n^{-1/\alpha}\left(\frac{1}{\alpha - 1} \sum_{k < \tau_n} f \bigg( \frac{X_k}{n} \bigg)
- n \int_0^1 f(x) \: dx\right),
\end{equation}
where $X=(X_k)_{k=0}^{\tau_n}$ is the Markov chain embedded into the block-counting process associated with the Beta$(2 - \alpha, \alpha)$-coalescent.
\begin{Theo}\cite[Theorem 2.1]{GKAW-KWS14}\label{Betafluc}
For all $f \in \mathscr F$
\begin{align*}   \mathcal J_n(f) = \int_0^\infty  f(m(r))m(r) d\mathcal L_{n,-r} + o_P(1),
\end{align*}
  where $m( r ) := \big( \frac{\alpha \Gamma(\alpha)}{r + \alpha \Gamma(\alpha)} \big)^{1/(\alpha - 1)}, \, r \ge 0$, and  $\mathcal L_n =\left(\mathcal L_{n,s}\right)_{s\in \mathbb R\}}$, is a compensated L\'evy process \index{process!Levy@L\'evy}\index{Levy@L\'evy process}with L\'evy measure $\beta_\alpha p^{-1-\alpha} dp$,  $\beta_\alpha := \frac 1{\Gamma(\alpha)\Gamma(2-\alpha)}$.
    \end{Theo}
    The statement of Theorem \ref{Betafluc} is an asymptotic version of the equation
\begin{align*}\label{sums}
n^{-1/\alpha} \left(\eta\sum_{k < \tau_n} f \bigg( \frac{X_k}{n} \bigg)  - n \sum_{k < \tau_n} f \bigg( \frac{X_k}{n} \bigg) \frac{\Delta_{k+1}}{n}\right) = \sum_{k < \tau_n} f \bigg( \frac{X_k}{n} \bigg)\frac{\Delta_{k+1}-\eta}{n^{1/\alpha}} , 
\end{align*}
where $\eta:= (\alpha-1)^{-1}$ and as above $\Delta_k=X_{k-1}-X_{k}$, $0\le k < \tau_n$. The proof of Theorem \ref{Betafluc} consists of showing that the  second and the third sum appearing in that equation may asymptotically be replaced by the  integrals $ \int_0^1 f(x) \: dx$ and $\int_0^\infty  f(m(r))m(r) d\mathcal L_{n,-r}$, respectively.
It is implicit in the statement of Theorem \ref{Betafluc} that the process $\mathcal L_n$ must be constructed on the same probability space as the \index{coalescent!n@$n$-}$n$-coalescent. We achieve this by constructing  
    the Poisson point measure (PPM) $\Psi_n$ of jumps of $\mathcal L_n$ from the PPM $\Phi$ with intensity measure $dt p^{-2} \Lambda(dt) = dt \,  \beta_\alpha p^{-1-\alpha} (1-p)^{\alpha-1} \mathbf 1_{p\le 1} dp$. This is done in  two steps. First, an independent Poissonian superposition is added to $\Phi$ to obtain a PPM $\Upsilon'$ on $\mathbb R \times \mathbb R_+$ with intensity measure $dt \times \beta_\alpha\ p^{-1-\alpha} dp$, see Fig.~\ref{GKAW-fig1}. (This as well as the next figure are borrowed from  \cite{GKAW-KWS14}.)

    \begin{figure}[h]
    \begin{center}
\includegraphics[width=10.5cm]{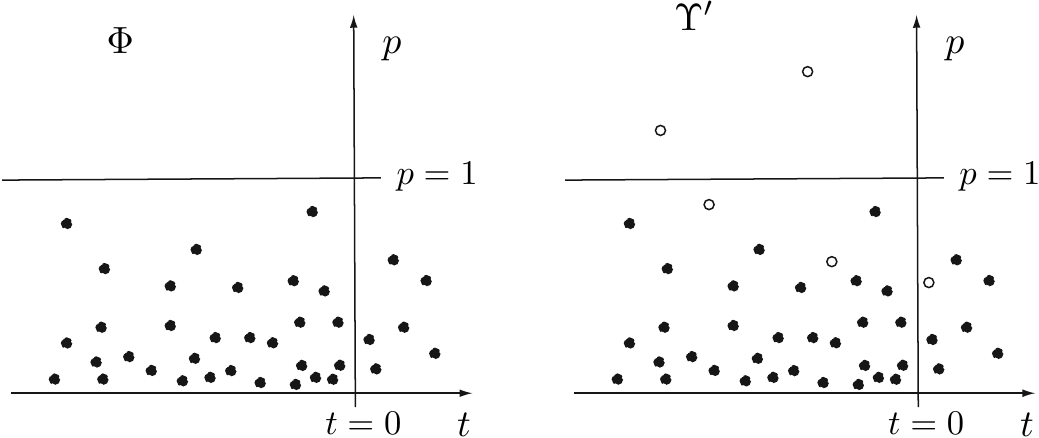}
\end{center}
\caption{\label{GKAW-fig1} 
left: PPM $\Phi$ with intensity   $dt \times \beta_\alpha p^{-1-\alpha}
(1-p)^{\alpha-1} \mathbf 1_{p\le 1} dp$  \newline right:  The PPM $\Upsilon'$ contains the points (some of them displayed and marked by~$\bullet$) from $\Phi$, plus additional points (marked by~$\circ$) that  make up for the difference between the intensities $p^{-1-\alpha} (1-p)^{\alpha-1}  1_{p\le 1}$ and $p^{-1-\alpha}$, $p>0$.}
\end{figure}

\begin{figure}[h]
\includegraphics[width=12cm]{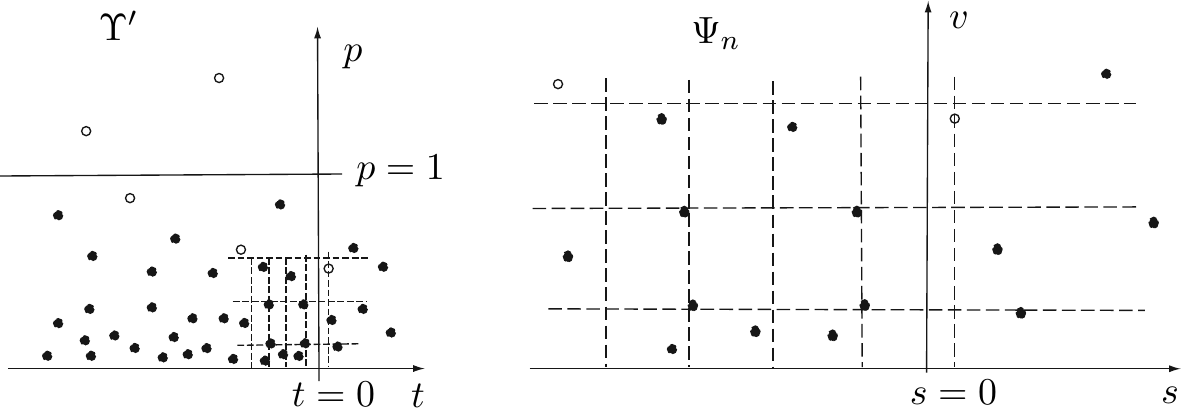} 
\caption{\label{GKAW-fig2}The point process $\Psi_n$ arises from $\Upsilon'$ through the transformation  $(t,p) \mapsto (s,v):=(n^{\alpha-1}t, n^{1-1/\alpha}p)$.}
\end{figure}

In a second step, as illustrated in Fig.~\ref{GKAW-fig2},  from $\Upsilon'$  we arrive at the PPM~$\Psi_n$ via the mapping 
\begin{equation}
\label{genevol}(t, p) \mapsto (s,v):= (n^{\alpha-1} t, n^{1-1/\alpha} p).
\end{equation}
It is easy to check that this transformation leaves the intensity measure, and hence also the distribution, of the PPM $\Phi$  invariant. Thus the distribution of the PPM $\Psi_n$, and hence also that of  the L\'evy process $\mathcal L_n$, does not depend on $n$.

The following statement is immediate from Theorem \ref{Betafluc}, and crucial for the proof of Corollary \ref{Cor2} on the fluctuations in the \index{evolving coalescent}\index{coalescent!evolving}evolving Beta-coalescent.

\begin{corollary} \cite[Corollary 2.2]{GKAW-KWS14}\label{mBetafluc} Let $f_1,\ldots, f_d \in \mathscr F$. Then
\[(\mathcal J_n(f_1),\ldots, \mathcal J_n(f_d)) \stackrel {\rm d}\to ( - I(f_1), \ldots, - I(f_d)),\]
where for $f\in \mathscr F$ we put $I(f) :=-\int_0^\infty  f(m(r))m(r) d\mathcal L_{-r}$, with $m$ being as in Theorem \ref{Betafluc} and $\mathcal L$ being a compensated \index{Levy@L\'evy process}\index{process!Levy@L\'evy}L\'evy process with $\mathcal L_0=0$ and  L\'evy measure $\beta_\alpha p^{-1-\alpha} dp$.
\end{corollary}

\begin{figure}[h]
\begin{center}
\includegraphics[width=12cm]{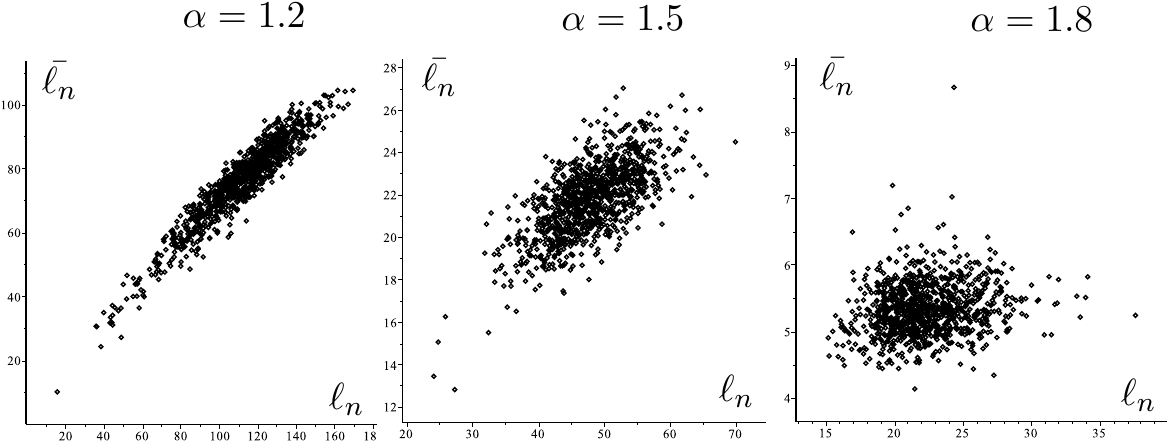}
\end{center}
\caption{\label{GKAW-fig4} Joint empirical distribution of the total length $\ell_n$ and the total external length $\bar{\ell}_n$ of an \index{Beta-coalescent}\index{coalescent!Beta-}$n$-Beta($2-\alpha, \alpha$)-coalescent for $n=1000$ and three exemplary values of $\alpha$ (one far below, one slightly below, one above the golden ration $g$).}  
\end{figure}

\begin{Exm}  \cite[Examples 2.6-2.8]{GKAW-KWS14} \label{Exam}
For $1 < \alpha < 2$  the fluctuations of the \index{branch length!total external}total external length, which first have been investigated in \cite{GKAW-DKeW14}, are captured by putting $f(x) = \alpha (\alpha - 1) (2 - \alpha) \Gamma(\alpha)$ in  Theorem   \ref{Betafluc}. They are of the order $O_P(n^{1-\alpha + 1/\alpha})$, which is $\gg 1$ for  $\alpha <g:= \frac 12 (1+\sqrt 5) $ (the ``golden ratio''), and of order $o_P(1)$  for $\alpha > g$. As discovered in \cite{GKAW-K12}, the fluctuations of the total length are   again $\gg 1$ for $\alpha < g$, but, in contrast, of order 1 for $\alpha > g$. In this latter parameter regime the total length's fluctuations  originate mainly from oscillations arising close to the root, which  do not show up in the total external length. For details see \cite{GKAW-K12}.

In the  case $\alpha <g$ the fluctuations of the total length  are captured by putting $f(x) :=$ \mbox{$\alpha (\alpha - 1) \Gamma(\alpha) x^{1 - \alpha}$} in Theorem   \ref{Betafluc}, and  Corollary \ref{mBetafluc} then describes the joint fluctuations of the total length and the total external length. Fig. \ref{GKAW-fig4} illustrates this by simulation results. 
\end{Exm}

So far we discussed fluctuation results on \index{Beta-coalescent}\index{coalescent!Beta-}Beta-coalescents under circumstances where randomness   arises mainly from the Markov chain $X$ embedded in  the block-counting process $N_n$. However, this scenario is inapplicable in important cases, most notably for the Kingman coalescent, where $X$ is a deterministic process. This requires specifically tailored methods. In the following instance this amounts to couple the quantities of interest to  certain Markov chains.

\begin{Theo}\cite[Main Theorem]{GKAW-DK15}\label{KingmanSFS}
Let $\bar \ell_{n,k}$ be the \index{tree length}\index{branch length!total}total length of all branches of order $k\ge 1$, that are those branches subtending  $k$ leaves. Then we have for the \index{Kingman coalescent}\index{coalescent!Kingman}Kingman coalescent as $n \to \infty$ 
\[ \sqrt{\frac{n}{4\log n}}\Big(\bar \ell_{n,k}-\frac 2k\Big) \stackrel {\rm d}\to  N(0,1)  \]
for each $k \ge 1$, and the  branch lengths of different order are asymptotically independent.
\end{Theo}

\noindent
A corollary asserting asymptotic Poissonian fluctuations for the site frequency spectrum of the Kingman coalescent is immediate, see \cite{GKAW-DK15}.

\section[Evolving $n$-coalescents]{Evolving $n$-coalescents and their limits}\label{evolvn}

\subsection{A Poisson construction of the evolving $n$-coalescent}
 \label{Poisfinite}
 As in the previous section we think of  $dt \,  p^{-2} \Lambda(dp)$ as the intensity measure of $p$-mergers, where $\Lambda$ is a finite measure on $(0,1]$. Let the population size $n\in \mathbb N $ be fixed; individuals will be labeled $1, \dots, n$. The following construction, which describes the evolving genealogy of the population driven by $\Lambda$, appears in~\cite{GKAW-KWS14}. It combines elements of the Poisson process constructions of the $\Lambda$-coalescent given in \cite{GKAW-P99}  and of the \index{evolving coalescent}\index{coalescent!evolving}\index{Bolthausen-Sznitman coalescent}\index{coalescent!Bolthausen-Sznitman}evolving Bolthausen-Sznitman coalescent in~\cite{GKAW-S12}.

Let $ \Upsilon=\Upsilon_n$ be a Poisson point process on $\mathbb R \times (0,1] \times [0,1]^n$ with intensity measure 
$$dt\, {p^{-2} \Lambda(dp) \, dv_1 \dots dv_n}.$$  Suppose $(t, p, v_1, \dots, v_n)$ is a point of $\Upsilon$.  If zero or one of the points $v_1, \dots, v_n$ is less than $p$, then no reproductive event occurs at time $t$.  However, if $k \geq 2$ of these points are less than $p$, so that $v_{i_1} < \dots < v_{i_k} < p$, then at time $t$, the individuals labeled $i_2, \dots, i_k$ all die, and the individual labeled $i_1$ gives birth to $k - 1$ new individuals who are assigned the labels $i_2, \dots, i_k$.  Seen backwards in time, this amounts to a coalescence of the lineages labeled $i_1, \dots, i_k$  at time $t$, and the rate of events that cause the lineages $i_1, \dots, i_k$ to coalesce is $\lambda_{n,k}.$  
\begin{figure}[h]
\begin{center}
\includegraphics[height=3cm]{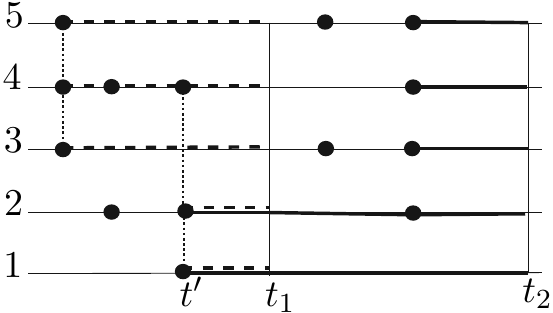}
\end{center}
\caption{\label{GKAW-fig3}The picture shows realisations of $\mathcal T_5(t_1)$ (dashed lines) and $\mathcal T_5(t_2)$ (solid lines), with the root of $\mathcal T_5(t_2)$ being at time $t'$. Each black dot marks a point $(t,i)$ with $(t,p,v_1,\ldots,v_5) \in \Upsilon$ such that $v_i < p$.  Only those coalescences are marked (by dotted lines) which affect $\mathcal T_5(t_1)$. The coordinates of the point $(t',p',v_1',\ldots,v_5') \in \Upsilon$, which leads to a merger both for $\mathcal T_5(t_1)$ and for $\mathcal T_5(t_2)$, obey $v_4' < \min(v_1',v_2') < \max(v_1',v_2') < p' < \min(v_3',v_5')$.}
\end{figure}
  For each $t \in \mathbb R$ one has a realization of the $n$-coalescent which describes the genealogy of the  $n$ individuals at time $t$.  We denote the corresponding coalescent tree (read off from the genealogy backwards from time  $t$)  by ${\mathcal T}_n(t)$ and call the process $({\mathcal T}_n(t), t \in \mathbb R)$ the \index{coalescent!evolving $n$-}\index{evolving coalescent}evolving $n$-coalescent.  Figure \ref{GKAW-fig3}, which is taken from \cite{GKAW-KWS14}, gives an illustration. Let us note that, unlike the \index{lookdown construction}lookdown construction (which will be discussed in Section~\ref{EvolvLambda}), this construction does not exhibit a pathwise consistency between different~$n$.
\subsection{Fluctuations in evolving \index{evolving coalescent}\index{coalescent!evolving $n$-}\index{Beta-coalescent}\index{coalescent!Beta-}Beta-coalescents}\label{evolvbeta}
Let us now specialize to the case $\Lambda={\rm Beta}(2-\alpha, \alpha)$, $1<\alpha < 2$.  For each  $s \in \mathbb R$ and $n \in \mathbb N$, we denote the Markov chain embedded into the block-counting process of  the coalescent tree $\mathcal T_n(n^{1-\alpha}s)$ by $(X_k^s)_{k=0}^{\tau_n^s}$. 
By shifting the origin of the scaled time to the time point~$s$ and re-centering the process $L_n$ at this new time origin (which does not affect its increments), we can apply Theorem~\ref{Betafluc} and conclude that
\[\frac{1}{\alpha - 1} \sum_{k < \tau_n^s} f \bigg( \frac{X_{k}^s}{n} \bigg)  = n \int_0^1 f(x) \: dx  + n^{1/\alpha} \int_{0}^\infty f(m(r))m(r)\; d\mathcal L_{n,s-r} + o_P(n^{1/\alpha}).
\]
Writing $\mathcal J_{n,s}(f) $ for the random variable \eqref{In} with $(X_k)=(X_k^0)$ replaced by $(X_k^s)$, we thus obtain
\begin{align*} 
\mathcal J_{n,s}(f) 
= \int_{0}^{\infty} f(m(r)) m(r)\: d\mathcal L_{n,s-r} + o_P(1).
\end{align*} 
The following result is now an immediate consequence of Corollary \ref{mBetafluc}.
\begin{Theo} \label{Cor2}
For $f \in \mathscr F$ and $s \in \mathbb R$, let $\mathcal J_{n,s}(f)$ be as in  \eqref{In}, but now evaluated at the \index{coalescent tree}coalescent tree $\mathcal T_n(n^{1-\alpha}s)$ instead of $\mathcal T_n(0)$. Then the sequence of stationary processes $(\mathcal J_{n,s}(f))_{-\infty < s < \infty}$,  $n \ge 1$,
converges as $n \to \infty$ in finite-dimensional distributions  to the moving average process
\[  \int_{0}^{\infty} f(m(r)) m(r)\: d\mathcal L_{s-r} , \hspace{.3in}  -\infty < s < \infty, \]
where $m$ and $\mathcal L$ are as in Corollary \ref{mBetafluc}.
\end{Theo}

\begin{remark} \label{remarkevtime}
As in Example \ref{Exam} this theorem applies to the total external branch length of \mbox{$n$-coalescents}, and to the tree length in the restricted parameter range $1<\alpha < g$.
Notably, the time scale $s$ is the \index{time scale!generation}{\em generation time scale} of the evolving \mbox{$n$-coales}\-cent: a single lineage is affected  by order of one reproduction event in one unit of this time scale, cf. Remark \ref{remarkgentime}. In contrast to this, $t$ is the \index{time scale!evolutionary}{\em evolutionary time scale}; in accordance with  Remark \ref{remarkgentime}, one time unit of the latter corresponds to $n^{\alpha-1}/(\alpha(\alpha-1)\Gamma(\alpha))$ units of the generation time scale. For simplicity, we have omitted the constant in the transformation \eqref{genevol}.
\end{remark}

Once again, the \index{coalescent!Kingman}\index{Kingman coalescent}\index{coalescent!evolving}\index{evolving coalescent}evolving Kingman coalescent does not fit into this framework. Here different functionals require varying techniques. For the case of  the \index{tree length}\index{branch length!total}tree length we refer to \cite{GKAW-PWW11}. The  total external branch length may be treated by means of 
Theorem \ref{KingmanSFS}. Thereby we obtain the following result.

\begin{Theo}\cite[Theorem 1]{GKAW-DK17}
Let $\bar \ell_n(s)$ be the \index{branch length!total external}total external branch length of the \index{coalescent tree}coalescent tree $\mathcal T_n(s/n)$, $s\in \mathbb R$. Then the sequence  of stationary processes \[\Big(\sqrt{\frac n{4\log n}}\Big(\bar \ell_n(s)-2\Big)\Big)_{-\infty < s < \infty}\ ,\] $n \ge 1$, converges as $n \to \infty$ in finite-dimensional distributions  to the stationary Gaussian process $(\bar \ell(s))_{-\infty < s < \infty}$ with mean zero and covariance function
\[\mathbb{COV}(\bar \ell(s), \bar \ell(0)) = \Big( \frac 2{2+s}\Big)^2 \ , \ s \in \mathbb R . \]
This process has a.s. continuous paths.
\end{Theo}

We point out that, in the corresponding theorem from \cite{GKAW-PWW11} on the total length, the time axis remains the (unscaled) evolutionary time $t$. This reflects the fact that, for the Kingman coalescent, the fluctuations of the total length come from contributions close to the root. In contrast, the fluctuations of the external length originate at the leaves and require the time scaling as given in the previous theorem.  On the other hand, this discrapancy of time scales disappears for Beta($2-\alpha,\alpha$)-coalescents with $1< \alpha < g$, see Remark \ref{remarkevtime}. Overall, this highlights the fact that fluctuations close to the root happen on the evolutionary time scale, whereas fluctuations close to the leaves fit to the \index{time scale!generation}generation time scale. The results of \cite{GKAW-K12,GKAW-PWW11} thus suggest the following 

\begin{conjecture} Let $ \ell_n(t)$ be the treelength of the \index{coalescent tree}coalescent tree $\left(\mathcal T_n(t)\right)_{t\in \mathbb R}$. Then for \index{Beta-coalescent}\index{coalescent!Beta-}\mbox{Beta{\em(}$2-\alpha,\alpha${\em)}-}coalescents with $ \frac 12 (1+\sqrt 5) = g<\alpha < 2$ the sequence of stationary processes 
\[\big( \ell_n(t) - cn^{2-\alpha} \big)_{-\infty < t< \infty} Ê\]
with $c=\Gamma(\alpha)\alpha (\alpha-1)/(2-\alpha)$ converge in finite-dimensional distribution to a stationary pure jump process.
\end{conjecture}

Together with Remark \ref{remarkevtime} this points to a discontinuity at $\alpha=g$ in the time scaling of the fluctuations of the total length.

\section[Evolving $\Lambda$-coalescents]{Evolving $\Lambda$-coalescents and the lookdown space}
\label{EvolvLambda}
 \subsection{From the lookdown graph to the lookdown space}\label{LDspace}
While the Poissonian construction described in Subsection \ref{Poisfinite} lacks the strong consistency required to connect different subpopulation sizes $n$, the so-called \index{lookdown construction}{\em lookdown construction} of Donnelly and Kurtz~\cite{GKAW-DK99} does have this property. This is based on a Poissonian construction which allows to cover the space $\mathbb R\times \mathbb N$ by {\em lineages} that are {\em branching} in forward time and {\em coalescing} in backward time. Here is a brief description based on ~\cite{GKAW-Gu18b}. While Gufler's work also includes the case of coalescents with \index{merger!simultaneous multiple}simultaneous multiple mergers (where \index{Lambda@$\Lambda$-coalescent}\index{coalescent!Lambda@$\Lambda$-}$\Lambda$ is replaced by a finite \index{Xi@$\Xi$-coalescent}\index{coalescent!Xi@$\Xi$-}measure $\Xi$ on $\{(x_1,x_2,\ldots):x_1\ge x_2\ge\cdots \ge 0, \sum x_k \le 1\}$) and allows for general initial states at time $t=0$, we restrict in this short expos\'e to the $\Lambda$-case, and to a stationary time situation, not least because this fits nicely to the framework and the concepts described in the previous parts.

As in Section \ref{asexBeta}, let $\Phi$ be the Poisson point measure on $\mathbb R \times (0,1]$ with intensity $dt \,p^{-2} \Lambda(dp)$. We first translate $\Phi$ to a random point configuration~$\eta_0$  on $\mathbb R \times \mathcal P(\mathbb N)$ as follows: For each point \mbox{$(t,p)\in \Phi$,} let $\left(B(i)\right)_{i\in \mathbb N}$ be a Bernoulli$(p)$-sequence, and let $\mathcal S_t:=\{i\in \mathbb N: B(i)=1\}$.  Thus, in words,  $\mathcal S_t$ is the set of ``levels'' which (seen in the backward time direction) take part in the $p$-merger and (seen in the forward time direction) are affected by the corresponding reproduction event.  We then have the {\em configuration of macroscopic reproduction events} $\eta_{\rm mac}:= \sum_{(t,p)\in \Phi} \delta_{(t,\mathcal S_t)}$.
  In addition, and independently, we have for every pair of levels $i<j$ a Poisson proint process $N^{ij}$ on $\mathbb R$ with rate $\Lambda(\{0\})$. The {\em configuration of binary reproduction events} is $\eta_{\rm bin}:= \sum_{i<j}\sum_{t\in N^{ij}} \delta_{(t,\{i,j\})}$.

We now describe the propagation of lineages at a (macroscopic or binary) reproduction event $(t,\mathcal S)\in \eta:= \eta_{\rm mac}\cup \eta_{\rm bin}$.
The individual that sits at level $i=\min \mathcal S$ at time $t-$ has offspring on all $j \in \mathcal S$ at time~$t$, and the lineages that at time $t-$ sit  at the second largest element of  $\mathcal S$ or above  are pushed up at time $t$  so that their order on the levels is preserved, they give way to the newly inserted lineages at $\mathcal S$,    and every level in $\mathbb N$ remains occupied by exactly one lineage.

For $s \le t$ we dentote by $A_s(t,i)$ the level of the ancestor at time $s$ of the individual $(t,i)$, and call $A(t,i)= (A_s(t,i))_{s\le t}$ the \index{ancestral lineage}{\em ancestral lineage} of $(t,i)$. In concordance with the above description, the maps $s\mapsto A_s(t,i)$ and $t\mapsto A_s(t,i)$ are c\`adl\`ag.
 For any two individuals $z_1=(t_1,i_1)$, $z_2=(t_2,i_2)$ there is a unique individual $(t,i)$ (with $i \le \min(i_1, i_2)$ and $t \le \min(t_1,t_2)$) in which the ancestral lineages $A(z_1)$ and $A(z_2)$ merge; we call  this individual the \index{most recent common ancestor}\index{dust-free}{\em most recent common ancestor} of the two individuals $z_1$ and $z_2$ and denote it by MRCA$(z_1, z_2)$. The genealogical distance of the two individuals is
\[\rho(z_1,z_2) := |t_1-t|+|t_2-t|.\]
The distance $\rho$ is a semi-metric on $\mathbb R \times \mathbb N$ (offspring individuals from the same parent have genealogical distance zero at the time of the reproduction event). We identify individuals with genealogical distance zero, and we take the metric completion. The resulting metric space $(Z, \rho)$ is called the \index{lookdown space}{\em lookdown space} \cite{GKAW-Gu18b} associated with $\eta$. In slight abuse of notation, we refer by $(t, i)\in \mathbb R \times \mathbb N$ also to the element of the metric space after the identification of elements with $\rho$-distance zero, in this sense we also assume $\mathbb R \times \mathbb N\subset Z$.

In the next subsection we will see that -- at least in the \index{dust-free}dust-free case -- the lookdown space can a.s.~be equipped with a family of probability measures $\mu_t$, $t \in \mathbb R$, which turns $(Z,\rho,\mu_t)$ into a family of  \index{metric measure space}metric measure spaces. There is no reasonable concept like  a ``space of all metric measure spaces''; however this kind of difficulty is resolved when passing to \index{metric measure space!isomorphy class}{\em isomorphy classes}. Two metric measure spaces $(X',r',\mu')$, $X'',r'', \mu'')$ are said to be {\em isomorphic} if there exists an isometry $\varphi$ from the closed support of $\mu'$ in $X'$ to the closed support of $\mu''$ in $X''$, with $\mu''=\varphi(\mu')$, the image of $\mu'$ under $\varphi$. The corresponding isomorphy class is denoted by $[[X', r', \mu']]$, and the space $\mathbb M$ of isomorphy classes of metric measure spaces is equipped with the \index{metric measure space!Gromov-Prohorov}Gromov-Prohorov metric $d_{\rm GP}$. This metric is complete and separable and, as shown in \cite{GKAW-GPW09}, induces the Gromov-weak topology, i.e. the weak topology on the distance matrix distributions when sampling independently from the probability measure~$\mu'$. 

There is yet another isomorphy concept between  metric measure spaces which is relevant for the concergence results of \cite{GKAW-Gu18b} that are reviewed in the next subsection. Two metric measure spaces $(X', r', \mu'), (X'', r'', \mu'')$ are called {\em strongly isomorphic} if they are measure-preserving isometric, that is, if there exists a surjective isometry $\varphi : X \to  X'$ with $\mu'' = \varphi(\mu')$. The {\em strong isomorphy class} is denoted by $[X',r',\mu']$. The set $\mathbf M$ of strong isomorphy classes of compact metric measure spaces can be equipped with a metric, the so-called \index{metric measure space!Gromov-Hausdorff-Prohorov}Gromov-Hausdorff-Prohorov metric $d_{\rm GHP}$. It turns $(\mathbf M, d_{\rm GHP})$ into a complete, separable metric space and induces the Gromov-Hausdorff-Prohorov topology on~$\mathbf M$ (see \cite{GKAW-EW06}). Because of the apperance of the Hausdorff distance in the definition of $d_{\rm GHP}$, and because of the requirement that the isometry works also on sets that are not charged by the sampling measures, the Gromov-Hausdorff-Prohorov metric controls more features than the Gromov-weak topology, as will become manifest in the following theorems.

 \subsection{Evolving $\Lambda$-coalescents \index{evolving coalescent}\index{coalescent!evolving}as continuum \index{tree-valued process}\index{process!tree-valued}tree-valued processes}\label{evolvcoal}
Let the Poisson random measure $\eta$ be defined from the finite measure $\Lambda$ as in the previous subsection, and let $(Z, \rho)$ be the lookdown space associated with $\eta$. We denote the Prohorov metric on the space of probability measures on $(Z,\rho)$ by $d_{\rm Proh}^Z$.
 For each $t\in \mathbb R$ and $n\in \mathbb N$, let the probability measure $\mu_t^n$ on $(Z,\rho)$ be the uniform measure on the first $n$ individuals at time $t$, that is,
 \[\mu_t^n = \frac1n \sum_{i=1}^n \delta_{(t,i)}.\]
 Let $\mu_t$ be the weak limit of $\mu_t^n$ as $n\to \infty$, provided the limit exists. We write 
\[\Theta_0 := \{t\in \mathbb R: \mbox{ there exists a } p \in (0,1] \mbox{ with } (t,p)\in \Phi\}.\]
for the set of times at which a  {\em macroscopic} reproduction event happens.
\begin{Theo}\label{treeprocnodust} \cite[Theorem 3.1, Propositions 4.1 and 4.7]{GKAW-Gu18b}
Assume the \index{Lambda@$\Lambda$-coalescent}\index{coalescent!Lambda@$\Lambda$-}$\Lambda$-coales\-cent is \index{dust-free}dust-free.Then there exists an event of probability one on which the following assertions hold:
\begin{itemize}
\item[i)]
For all $T \in \mathbb R$,
\[\lim_{n\to \infty} \sup_{t\in[0,T]}d_{\rm Proh}^Z(\mu_t^n,\mu_t) = 0.\]
\item[ii)] The map $t \mapsto \mu_t$ is c\`adl\`ag in the weak topology on the space of probability measures on $(Z,\rho)$, and  the set $\Theta_0$ is the set of jump times.
\item[iii)] For all $t\in \Theta_0$, the measure $\mu_t$ contains atoms, and the left limit $\mu_{t-}$ is non-atomic.
\end{itemize}
Moreover, a Feller-continuous strong Markov process with values in $\mathbb  M$ is given a.s. by $([[Z, \rho, \mu_t]], t\in \mathbb R)$. This process has c\`adl\`ag paths in the \index{metric measure space!Gromov-Prohorov}Gromov-Prohorov metric and $\Theta_0$ is the set of its jump times.
\end{Theo}
For $(s,i) \in \mathbb R \times \mathbb N$ and $t \ge s$  let  
\[D_t(s,i) := \min\{j\in \mathbb N:A_s(t,j)=i\}.\]
Then $D(s,i) := (D_t(s,i))_{t\ge s}$ is the leftmost {\em line of descent} of $(s,i)$. This is the longest-living among all \index{line of descent}lines of descent of $(s,i)$; it corresponds to the fixation line in \cite{GKAW-H15} and to the forward level process in \cite{GKAW-PW06} (where  the Kingman case $\Lambda = \delta_0$ had been studied with a focus on the birth and fixation times of the \index{most recent common ancestor}consecutive MRCA's of the entire population). The map $t\mapsto D_t(s,i)$ is non-decreasing. Let 
\[\tau_{s,i}=\inf\{t\in [s,\infty): D_t(s,i) = \infty\},\]
in other words, $\tau_{s,i}$ is the extinction time of the offspring of the individual $(s, i)$.
Then the set  of times at which such ``old families'' become extinct is given by 
\[\Theta^{\rm ext} := \{\tau_{s,i}: s\in \mathbb R,\, i\in \mathbb N\}.\]
In the Kingman case, the set $\Theta^{\rm ext}$ is described in \cite[Proposition~1]{GKAW-DKW14}  as a superposition of Poisson point processes.

For $t\in \mathbb R$, let $Z_t$ be the closure of the set $\{t\}\times \mathbb N$, seen as a subspace of the complete space $(Z,\rho)$. If the $\Lambda$-coalescent comes down from infinity, then there exists an event of probability one on which all subsets $Z_t$, $t \in \mathbb R$,  are compact. Indeed, $|\{A_s(t,j): j\in \mathbb N\}|<\infty$ for all $s< t$ implies, by definition of the metric $\rho$,  that the complete subspace $Z_t$ is totally bounded in $Z$.
\begin{Theo}\label{treeprocCDI} \cite[Theorem 3.5, Propositions 4.2 and 4.8]{GKAW-Gu18b}
Assume that the \index{coalescent!Lambda@$\Lambda$-}\index{Lambda@$\Lambda$-coalescent}$\Lambda$-coales\-cent \index{coming down from infinity}comes down from infinity. Then the following assertions hold on an event of probability one:
\begin{itemize}
\item[i)]
For each $t \in \mathbb R$, the compact set $Z_t$ is the closed support of $\mu_t$.
\item[ii)] The map $t \mapsto Z_t$ is c\`adl\`ag for the Hausdorff distance on the set of closed subsets of $(Z,\rho)$. The set $\Theta^{\rm ext}$ is the set of jump times. For each $t\in \Theta^{\rm ext}$, the set $Z_t$ and the left limit $Z_{t-}$  are not isometric.
\end{itemize}
Moreover, a Feller-continuous strong Markov process with values in $\mathbf  M$ is given a.s. by $([Z, \rho, \mu_t], t\in \mathbb R)$. It has c\`adl\`ag paths in the \index{metric measure space!Gromov-Hausdorff-Prohorov}Gromov-Hausdorff-Prohorov topology and $\Theta_0 \cup \Theta^{\rm ext}$ is the set of its jump times.
\end{Theo}
Next we turn to the case with dust. For $(t,i)\in \mathbb R\times \mathbb N$ let us denote by $z(t,i)$ the most recent ancestor of $(t,i)$ that participates {\em together with other individuals at the same time} in a merging event. This individual is of the form $z(s,j)$ for some $s\le t$. We define $v(t,i):= \rho((t,i), z(t,i))= t-s$, that is the length of the external branch in the infinite coalescent back from the leaf $(t,i)$. We endow $Z\times \mathbb R_+$ with the product metric $d^{Z\times \mathbb R_+}((z,v) ,(z',v'))=\rho(z,z')\vee |v-v'|$, and denote the Prohorov metric on the space of probability measures on $Z\times \mathbb R_+$ by $d_{\rm Proh}^{Z\times \mathbb R_+}$. For each $t\in \mathbb R$ and $n \in \mathbb N$, we define a probability measure  on $Z\times \mathbb R_+$   by 
\[m_t^n = \frac 1n \sum_{i=1}^n \delta_{(z(t,i),v(t,i))}.\]
Let $m_t$ denote the weak limit of $m_t^n$ as $n\to \infty$, provided that limit exists. The following theorem guarentees that the family $(m_t)$ is well-defined a.s., thus $(Z,\rho, m_t), \, t\in \mathbb R$, constitutes a family of so-called \index{metric measure space!marked}{\em marked metric measure spaces} with mark space $\mathbb R_+$.  Let $\widehat {\mathbb M}$ be the set of isomorphy classes of such marked metric measure spaces. $\widehat {\mathbb M}$ is endowed  with the \index{metric measure space!marked Gromov-Prohorov}{\em marked Gromov-Prohorov metric} $d_{\rm mGP}$, which makes $(\widehat {\mathbb M}, d_{\rm mGP})$ a complete, separable metric space (\cite{GKAW-DGP11}).  Again, the isomorphy class of $(Z,\rho, m_t)$ will be denoted by $[[Z,\rho, m_t]]$.
\begin{Theo}\label{treeprocdust} \cite[Theorem 3.10, Propositions 4.12 and 4.15]{GKAW-Gu18b}
Assume that the \index{Lambda@$\Lambda$-coalescent}\index{coalescent!Lambda@$\Lambda$-}$\Lambda$-coales\-cent has a \index{dust component}dust component. Then the following assertions hold on an event of probability one:
\begin{itemize}
\item[i)]
For all $T \in \mathbb R$,
\[\lim_{n\to \infty} \sup_{t\in[0,T]}d_{\rm Proh}^{Z\times \mathbb R_+}(m_t^n,m_t) = 0.\]
\item[ii)] The map $t \mapsto m_t$ is c\`adl\`ag in the weak topology on the space of probability measures on $Z\times \mathbb R_+$. The set $\Theta_0$ is the set of jump times. 
\item[iii)] For each $t\in (0,\infty)$, the left limit $m_{t-}$ satisfies $m_t(Z\times\{0\})=0$. For each $t\in \Theta_0$ one has $m_t(Z\times\{0\}) > 0$. 
\end{itemize}
Moreover,  a Feller-continuous strong Markov process with values in  $\widehat {\mathbb M}$ is given a.s. by $([[Z,\rho, m_t]], t\in \mathbb R_+)$. This process has a.s.~c\` adl\`ag paths in the \index{metric measure space!marked Gromov-Prohorov}marked Gromov-Prohorov metric and $\Theta_0$ is the set of its jump times.
\end{Theo}

The  constructions discussed above refer to {\em neutral} situations.  
We point out, however, that a pathwise approach which builds on the \index{lookdown space!with competition and selection}lookdown space turns out to be successful also in the presence of  selection and competition, and for fluctuating population sizes; this is work in progress joint by A. Blancas, S. Gufler, S. Kliem, V. C. Tran and the second author. 
We emphasise that the pathwise approach is just {\em one} way to construct \index{process!tree-valued}\index{tree-valued process}tree-valued processes, see Hutzenthaler\&Pfaffelhuber \cite{GKAW-HP20}  and Sturm\&Winter \cite{GKAW-SW20}  (and references therein), where tree-valued processes are obtained as solutions of well-posed martingale problems. The former  investigates tree-valued processes under fluctuating selection, and the latter studies multitype branching models with state-dependent mutation and competition. Also via well-posed martingales,  A. Winter \cite{GKAW-Wi20} treats tree-valued processes (and functionals thereof) from the viewpoint of algebraic measure trees.
For a probabilistic analysis of the so-called common ancestor type distribution in $\Lambda$-evolutions with selection and mutation we  refer to~\cite{GKAW-BLW16}.

\enlargethispage*{1cm}

\paragraph{\bf Acknowledgement.} We thank two anonymous referees for a careful reading and  helpful comments.

\end{document}